\documentclass[12pt]{elsarticle}
\usepackage{amsfonts}
\usepackage{amssymb}
\usepackage{mathrsfs}
\usepackage{amsmath}
\usepackage{amsthm}
\usepackage{hyperref}
\def\hD{\hat D}
\def\hP{\hat P}
\def\hM{\hat M}

\def\hX{\hat X}

\def\cL{\mathcal L}
\def\cP{\mathcal P}
\def\cG{\mathcal G}
\def\cD{\mathcal D}
\def\cF{\mathcal F}
\def\Z{\mathbb{Z}}
\def\N{\mathbb{N}}

\def\tr{^{\mathsf{T}}}
\def\eop{\unskip\nobreak\hfil\penalty50\hskip2em\hbox{}\nobreak
\hfill\mbox{$\Box $}\par}
\newtheorem{thm}{Theorem}[section]
\newtheorem{prop}{Proposition}[section]

\begin{document}
\begin{frontmatter}
	\title{ A Linear Algebra approach to monomiality and operational methods}
	
\author{Luis Verde-Star}
 \address{
Department of Mathematics, Universidad Aut\'onoma Metropolitana, Iztapalapa,
Apartado 55-534, Mexico City 09340,
 Mexico }
\ead{verde@xanum.uam.mx}

 \begin{abstract}
We use linear algebraic methods to obtain general results about linear operators on a space of polynomials that we apply to the operators associated with a polynomial sequence by the monomiality property. We show that all such operators are differential operators with polynomial coefficients of finite of infinite order. We consider the monomiality operators associated with several classes of polynomial sequences, such as Appell and Sheffer, and also orthogonal polynomial sequences that include the Meixner, Krawtchouk, Laguerre, Meixner-Pollaczek, and Hermite families.



 \end{abstract}
\end{frontmatter}

\section{Introduction}

Polynomials in one variable are among the simplest and most useful functions and
 sequences of polynomials have been studied since the times of Euler and Bernoulli. In particular, sequences 
 $\{p_n(t)\}_{n\ge 0}$ in one real or complex variable such that $p_n$ has degree $n$, for $n \ge 0$, which are called {\em polynomial sequences}, play a central role in several areas of mathematics and its applications. Properties of polynomial sequences related to combinatorics, numerical analysis, approximation theory, probability, coding theory, etc. have been extensively studied using diverse approaches. Around 1996 G. Dattoli and some collaborators \cite{DCS97} introduced operational methods based on the concept of monomiality, or quasi-monomiality. A polynomial sequence $\{p_n(t)\}_{n\ge 0}$ has the {\em monomiality property} if there exist linear operators $\hP$ and $\hM$ that satisfy 
\begin{eqnarray}\label{defMono}
	\hP p_n(t) &=& n p_{n-1}(t), \qquad n\ge 0, \cr
	\hM p_n(t) &=& p_{n+1}(t), \qquad n \ge 0.
\end{eqnarray}
The concept of monomiality has its origin in the poweroids of Steffensen, which appear in a natural way in the theory of finite differences. See the recent book \cite[p. 62--64]{LiccDat} for more information about the genesis of the monomiality principle.

A large number of papers that use monomiality to study diverse types of polynomial sequences have been published during the last 25 years. The methods of the umbral calculus are often used in combination with monomiality. See \cite{LiccDat} and the references therein. In some of those papers several generalizations of Hermite, Laguerre and other orthogonal polynomial families have been obtained. 

In this paper, we use linear algebraic methods to obtain general results about linear operators on the vector spaces of polynomials and formal power series. Then we apply the general results to study the operators related with the monomiality principle. We obtain our results working with the elements of an algebra $\cL$ of infinite matrices that contains all the lower triangular infinite matrices. Each element of $\cL$ can be interpreted as a sequence of polynomials and also as a linear operator on the space of polynomials. Certain particular and simple elements of $\cL$ represent the operators of differentiation and multiplication by the independent variable on the vector spaces of polynomials and formal power series and may be considered as generators of the algebra $\cL$.

Some results related with monomiality can be easily obtained in our algebraic setting. For example, the proof of the fact that every polynomial sequence has the monomiality property is almost trivial. We also show that all the lower triangular matrices in $\cL$ represent differential operators of infinite order with polynomial coefficients, and this also holds for all the operators $\hP$ and $\hM$ associated with a polynomial sequence.

Our main objective in this paper is to show that our linear algebraic approach can produce relevant results in the theory of monomiality and operational methods. For this purpose we present examples where we show how our methods are applied to several classes of polynomial sequences, such as those of Appell and Sheffer, and also some orthogonal polynomial sequences.

In our previous papers \cite{Opm} and \cite{Uni} we used our approach to study hypergeometric and basic hypergeometric orthogonal polynomial sequences, and in \cite{Lin} we showed how to obtain linearization and connection coefficients for general polynomial sequences.

This paper is organized as follows. In Section 2 we introduce the algebra $\cL$ of generalized Hessenberg matrices and present the basic definitions and properties of the matrices that we use in the rest of the paper. In Section 3 we show that every lower triangular matrix can be written as a matrix that represents a differential operator of infinite order with polynomial coefficients and how this results applies to the monomiality operators of arbitrary polynomial sequences. The class of Appell polynomial sequences is considered in Section 4.
In Section 5 we consider polynomial sequences of binomial type, which correspond to modified composition matrices. Sheffer sequences are studied in Section 6. 
In Section 7 we deal with a class of orthogonal polynomial sequences that includes the Meixner, Krawtchouk, Meixner-Pollaczek, Laguerre, and Hermite families of orthogonal polynomials. Finally, in Section 8 we present a simple way to extend the theory by replacing the differentiation operator with generalized derivatives, such as the Jackson $q$-derivative and Ward derivatives.

\section{ The algebra of generalized lower Hessenberg matrices}\label{algcL}

In this section we define an algebraic setting in which we can study polynomial sequences, sequences of formal power series, and linear operators acting on spaces of polynomials or formal power series. Here we define the algebra $\cL$ of generalized lower Hessenberg matrices and present some of their basic properties that will be used in the rest of the paper. 
 For additional information about the algebra $\cL$ see \cite{Opm} and \cite{PSGH}.

An infinite matrix $A=[a_{j,k}]$, where the indices run over the non-negative integers and the entries are complex numbers is a {\sl lower generalized Hessenberg matrix} if there exists an integer $m$ such that $a_{j,k}=0$ whenever
 $j-k < m$. We denote by $\cL$ the set of all such matrices.
 
 We say that the entry $a_{j,k}$ of $A$ lies in the $n$-th diagonal of $A$ if $j-k=n$. If $m > n$ then the $m$-th diagonal lies below (to the left of) the $n$-th diagonal.
 A nonzero element of $\cL$ is a {\sl diagonal matrix} if all of its nonzero elements lie in a single diagonal. 
 If $A$ is a nonzero element of $\cL$ and $m$ is the minimum integer such that $A$ has at least one nonzero entry in the $m$-th diagonal, then we say that $A$ has {\sl index} $m$ and write ind$(A)=m$. The index of the zero matrix is infinity, by definition.

 It is clear that $\cL$ is a complex vector space with the usual addition of matrices and multiplication by scalars. It is also closed under matrix multiplication.
If $A$ and $B$ are in $\cL$, with ind$(A)=m$ and ind$(B)=n$, then the product $C=AB$ is a well defined element of $\cL$ and
\begin{equation}\label{eq:mult}
	c_{i,k}= \sum_{j=k+n}^{i-m} a_{i,j} b_{j,k}, \qquad i-k \ge m+n.
\end{equation} 
Note that ind$(AB)\ge m +n$ and that the multiplication in $\cL$ involves only finite sums.

 A sufficient, but not necessary, condition for $A$ to have a two-sided inverse is that ind$(A)=0$ and $a_{k,k}\ne 0$ for $k \ge 0$. We denote by $\cG$ the set of all matrices that satisfy such condition. It is clear that $\cG$ is a group under matrix multiplication. The unit is the identity matrix $I$, whose entries on the 0-th diagonal are equal to 1 and all other entries are zero.

We define next some particular elements of $\cL$ that will be used often in the rest of the paper. Let $X$ denote the diagonal matrix of index -1 with $X_{j,j+1}=1$ for $j \ge 0$, and denote by $\hat X$ the transpose of $X$. Note that $\hat X$ is diagonal of index 1, $X \hat X= I$ and $\hat X X= J_0$, where $J_0$ is the diagonal matrix of index zero that has its entry in the position $(0,0)$ equal to zero and its entries in positions $(j,j)$ equal to 1 for $j \ge 1$.
Therefore, $\hat X$ is a right-inverse for $X$, but it is not a left-inverse.

We say that a matrix $A \in \cL$ of index $m$ is {\sl monic} if all the entries of $A$ in the diagonal of index $m$ are equal to 1. Note that a monic matrix of index -1 is a unit lower Hessenberg matrix.

If $m$ is a positive integer then $X^m$ is the diagonal matrix of index $-m$ with all its entries in the $(-m)$-th diagonal equal to 1. Analogously, $\hX^m$ is diagonal of index $m$ and all its entries in the $m$-th diagonal are equal to 1.

We write $H=\hbox{Diag}(h_0,h_1,h_2,\ldots)$ to describe a diagonal matrix $H$ of index zero with $H_{k,k}=h_k$ for $k \ge 0$. 
Let $\cD$ be the group of invertible diagonal matrices of index zero. Let us define $N=\hbox{Diag}(1,2,3,\ldots)$ and $F=\hbox{Diag}(0!,1!,2!,\ldots)$. It is clear that $F$ and $N$ are elements of $\cD$.

Define $D= F \hX F^{-1}$. Note that $D$ is diagonal of index one, $D_{k+1,k}= k+1$ for $k \ge 0$, $X D=N$, and $D X=\hbox{Diag}(0,1,2,3,\ldots)$. Let $\hD$ be the transpose of $D$. Then $\hD$ is diagonal of index minus one. The notation $\hX$ and $\hD$ to denote  the transposes of $X$ and $D$ will only be used for these two matrices.  

We also have the commutation relation
\begin{equation}\label{eq:commutXD}
 X D - D X = I,
\end{equation}
where $I$ denotes the identity matrix $I=\hbox{Diag}(1,1,1,\ldots)$.

From  $ X D - D X= I$ it is easy to obtain the Pincherle differentiation formulas
\begin{equation}\label{eq:PinchDj}
	X D^j - D^j X = j D^{j-1}, \qquad j\ge 0, 
\end{equation}
and
\begin{equation}\label{eq:PinchXj}
	 X^j D - D X^j = j X^{j-1}, \qquad j\ge 0, 
	\end{equation}
and using induction we obtain 
\begin{equation}\label{eq:XkDj}
	X^k \dfrac{D^j}{j!}=\sum_{i= 0}^{\min(j,k)} \binom{k}{i} \dfrac{D^{j-i}}{(j-i)!} X^{k-i}, \qquad k, j \in \N.
\end{equation}

If $B$ is a diagonal matrix of index $m >0$ then $B^k$ is diagonal of index $k m \ge 0$ for $k \ge 0$ and therefore $\sum_{k\ge 0} a_k B^k$ is a well defined element of $\cL$ for any sequence of numbers $\{a_k\}_{k\ge 0}$. 

If $C$ is a diagonal matrix of index $m<0$ then $\sum_{k=0}^n a_k B^k$ is a banded element of $\cL$ of index at least equal to $n m$ for any finite sequence of numbers $a_0,a_1,\ldots, a_n$. 

The following almost trivial Proposition will be used later.

\begin{prop}\label{basic} 
Let $A$ be an element of the group $\cG$ and let $U$ be an element of $\cL$. Then we have

	(i) There exists a unique $U_r$ in $\cL$ such that $ A U_r= U A$.

	(ii) There exists a unique $U_\ell$ in $\cL$ such that $ U_{\ell} A = A U$. 
\end{prop}
{\it Proof.} Let $U_r=A^{-1} U A$ and $U_{\ell}= A U A^{-1}$. Then we have $ A U_r = A A^{-1} U A=UA$ and $ U_{\ell} A= A U A^{-1} A= A U$. \eop

 \subsection{Sequences of polynomials and formal power series} 

A matrix $A=[a_{k,j}]$ of index $m $ determines a sequence of polynomials $u_k(t)$ defined by 
\begin{equation}\label{eq:rowpoly}
u_k(t) = \sum_{j=0}^{k-m} a_{k,j} t^j, \qquad k \ge 0.
\end{equation}
That is, the entries in the $k$-th row of $A$ are the coefficients of $u_k$ with respect to the basis of monomials $\{t^k\}_{k\ge 0}$. 
If $m> 0$ then $u_j(t) =0$ for $ 0 \le j \le m-1$, and if $A$ is monic then each nonzero $u_k$ is monic and has degree $k-m$.
If $A$ is in the group $\cG$ then $m=0$ and $u_k(t)$ has degree $k$, for $k\ge 0$. That is, $\{u_k(t)\}_{k\ge 0}$ is a {\em polynomial sequence}.

Every matrix $B=[b_{k,j}]$ in $\cL$ determines a linear operator $\hat{B}$ on the vector space $\cP$ of polynomials defined by 
\begin{equation}\label{eq:maponP}
	\hat{B} t^k=\sum_{j\ge 0} b_{k,j} t^j, \qquad k \ge 0,
\end{equation}
that is, the image of $t^k$ is the polynomial associated with the $k$-th row of $B$. 
Therefore, if $A$ is in $\cL$ then multiplication on the right by $B$ can be interpreted as applying the map $\hat{B} $ to each of the polynomials associated with the rows of $A$. 

Each column of a matrix $A\in \cL$ determines a formal power series $f_j(x)$ defined by
\begin{equation}\label{eq:colseries}
f_j(x)= \sum_{k=0}^\infty a_{k,j} x^k, \qquad j \ge 0.
\end{equation}
If $A$ is in $\cG$ then $f_j(x)$ has lower degree equal to $j$, for $j \ge 0$. 

Every matrix $B=[b_{k,j}]$ in $\cL$ determines a linear operator $\check{B}$ on the vector space $\cF$ of formal power series $f(x)=\sum_{j=0}^\infty c_j x^j$.
The operator $\check{B}$ is defined by
\begin{equation}\label{eq:maponF}
	\check{B} x^j = \sum_{k\ge 0} b_{k,j} x^k, \qquad j \ge 0. 
\end{equation}
If a series $f(x)=\sum_{k\ge 0} c_k x^k$ is represented by the column vector $c= [c_0,c_1,c_2,\ldots]\tr $ then 
\begin{equation}
	\check{B} f(x) = \sum_{k\ge 0} e_k x^k, 
\end{equation}
where
$$	[e_0,e_1, e_2, \ldots ]\tr = B c, $$
that is
$$e_k= \sum_{j\ge 0} b_{k,j} c_j, \qquad k \ge 0.$$
These sums are finite since only a finite number of entries in the $k$-th row of $B$ can be nonzero. 
Therefore, if $A$ is in $\cL$ then multiplication on the left by $B$ can be interpreted as applying the map $\check{B}$ to each of the formal power series associated with the columns of $A$. 

\subsection{The basic shifts and differentiation operators}

Let $A=[a_{k,j}]$ be in $\cL$, let $u_k(t)$ be the polynomials associated with the rows of $A$, and let $f_j(x)$ be the formal power series associated with the columns of $A$. The matrices $X$ and $\hX$ defined previously act as shifts on the entries of $A$ as follows. 

(i) The rows of $A X$ are the rows of $A$ shifted one place to the right, with initial entry equal to zero.

(ii) The rows of $A \hX$ are the rows of $A$ shifted one place to the left.

(iii) The columns of $ X A$ are the columns of $A$ shifted one place upwards.

(iv) The columns of $\hX A$ are the columns of $A$ shifted one place downwards, with initial entry equal to zero.

From (i) we see that the polynomial corresponding to the $k$-th row of $AX$ equals $t u_k(t)$.
From (ii) it follows that the polynomial of the $k$-th row of $A \hX$ is equal to $(u_k(t)-u_k(0))/t$.
From (iii) we see that the series of the $j$-th column of $XA$ is equal to $(f_j(x) - f_j(0))/x$, and 
 (iv) says that the series of the $j$-th column of $\hX A$ is equal to $ x f_j(x)$.

The matrices $D$ and $\hD$, defined previously, act on $A$ as follows.
 For $k\ge 0$ and $j\ge 0$ we have 

(1) $(A D)_{k,j}= (j+1) a_{k,j+1}$. 

(2) $(A \hD)_{k,j}= j a_{k,j-1}$. 

(3) $(D A)_{k,j}= k a_{k-1,j}$. 

(4) $(\hD A)_{k,j}= (k+1) a_{k+1,j}$.

From (1) it follows that the polynomial of the $k$-th row of $AD$ is equal to $u_k^\prime(t)$.
From (2) we see that that the polynomial of the $k$-th row of $A \hD$ is equal to $t (t u_k(t))^\prime$.
From (3) we see that the series of the $j$-th column of $DA$ is equal to $ x (x f_j(x))^\prime$, and (4) gives us that the series of the $j$-th column of $\hD A$ is $f_j^\prime(x)$.

We can also describe the action of the operators considered above by their effects on the sequence of polynomials $\{u_k(t)\}_{k\ge 0}$ and the sequence of formal power series $\{f_j(x)\}_{j\ge 0}$ associated with $A$. From property (iii) we see that multiplying $A$ by $X$ on the left sends $u_k(t)$ to $u_{k+1}(t)$, for $k\ge 0$, and property (3) says that multiplying $A$ on the left by $D$ sends $u_k(t)$ to $ k u_{k-1}(t)$, for $k\ge 0$. In a similar way we can see from (2) that multiplying $A$ by $\hD$ on the right sends $f_j(x)$ to $ j f_{j-1}(x)$, for $j\ge 0$, and from (ii) that multiplying $A$ by $\hX$ on the right sends $f_j(x)$ to $f_{j+1}(x)$, for $j \ge 0$. 

\section{Differential operators on the space of polynomials}

In this Section we show that every linear operator on the space $\cP$ of polynomials that is represented by a lower triangular matrix is a differential operator of infinite order with polynomial coefficients. Such operators send $t^n$ to a polynomial of degree at most $n$, for $n \ge 0$. We also present a representation of the dual basis of a polynomial sequence associated with a matrix $A$ that uses the matrix $A X A^{-1}$, and some properties of the matrices $M$ and $P$ corresponding to a general element $A$ of the group $\cG$. 

\begin{thm}\label{diffOper}
	Let $A=[a_{i,j}]$ be a matrix of index zero. Then there exists a sequence of polynomials $\{p_k(t)\}_{k\ge 0}$ such $p_k(t)$ has degree at most equal to $k$, for $k\ge 0$, and 
	\begin{equation}
		A= \sum_{k\ge 0} \dfrac{ D^k}{k!} p_k(X).
	\end{equation}
\end{thm}
{\it Proof.} Let $r$ and $k$ be integers such that $k\ge r \ge 0$. The matrix $\frac{1}{k!} D^k X^{k-r} $ is diagonal of index $r$ and it is easy to see that its $(i, i-r)$ entry equals $\binom{i}{k}$ for $i \ge 0$. Define 
\begin{equation}\label{eq:dkr}
d(k,r)= \sum_{j=r}^k \binom{k}{j} (-1)^{k-j} a_{j,j-r}.
\end{equation}
The matrix 
$$E_r= \sum_{k\ge 0} d(k,r) \dfrac{D^k}{k!} X^{k-r}, \qquad r\ge 0, $$
is a well-defined diagonal matrix of index $r$. We show next that $E_r$ coincides with the diagonal of index $r$ of $A$. Let $i$ be a non-negative integer. 
Then we have
\begin{eqnarray*}
	(E_r)_{i,i-r} &=& \sum_{k\ge 0} d(k,r) \binom{i}{k} \cr
	 &=& \sum_{k\ge 0} \sum_{j=r}^k \binom{k}{j} \binom{i}{k} (-1)^{k-j} a_{j,j-r } \cr
	 &=& \sum_{j=r}^i a_{j,j-r} \sum_{k=j}^i \binom{i}{j} \binom{i-j}{k-j} (-1)^{k-j}\cr
	 &=& a_{i,i-r}.
\end{eqnarray*}
The last equality follows from
$$\sum_{k=j}^i \binom{i-j}{k-j} (-1)^{k-j}= (1-1)^{i-j}=\delta_{i,j}.$$
Therefore $E_r$ coincides with the diagonal of index $r$ of $A$ and then we have
\begin{equation}\label{eq:sumEr}
A= \sum_{r\ge 0} E_r=\sum_{r \ge 0} \sum_{k\ge 0} d(k,r) \dfrac{1}{k!} D^k X^{k-r}.
\end{equation}

Define 
\begin{equation}\label{eq:pkt}
p_k(t) = \sum_{r=0}^k d(k,r) t^{k-r}, \qquad k \ge 0.
\end{equation} 
Changing the order of the sums in \ref{eq:sumEr} we obtain
$$A= \sum_{k\ge 0} \dfrac{D^k}{k!} p_k(X).$$ \eop

Let $\{u_k(t)\}_{k\ge 0}$ be the sequence of polynomials associated with $A$, that is,
$$u_k(t)=\sum_{j=0}^k a_{k,j} t^j,\qquad k \ge 0.$$
Using the definitions \ref{eq:dkr} and \ref{eq:pkt} it is easy to obtain the inverse relations
\begin{eqnarray}\label{eq:binomTrans}
	p_k(t)&=& \sum_{j=0}^k \binom{k}{j} u_j(t) (-t)^{k-j}, \cr
        u_k(t)&=& \sum_{j=0}^k \binom{k}{j} p_j(t) t^{k-j}. 
\end{eqnarray}
These relations can be expressed in terms of generating functions as follows:
\begin{equation}\label{genFct}
\sum_{k=0}^\infty p_k(t) \dfrac{x^k}{k!}=	\exp(-t x) 	\sum_{k=0}^\infty u_k(t) \dfrac{x^k}{k!}.
\end{equation}

For $i\ge 0$ the sequence of polynomials associated with $A \frac{D^i}{i!}$ is $\left\{\frac{u_j^{(i)}(t)}{i!}\right\}_{j\ge 0}$. By the previous Theorem and \ref{eq:binomTrans} we have
\begin{equation}\label{eq:Aski}
	A \dfrac{D^i}{i!}= \sum_{k=0}^\infty \dfrac{D^k}{k!} s_{k,i}(X), \qquad i\ge 0,
\end{equation}
where
\begin{equation}\label{eq:ski}
s_{k,i}(t)= \sum_{j=i}^k \binom{k}{j} \dfrac{u_j^{(i)}(t)}{i!} (-t)^{k-j}.
\end{equation}
Note that $s_{k,i}(t)=0$ if $i>k$.

Using \ref{eq:binomTrans} we can write $s_{k,i}(t)$ in terms of the polynomials $p_j(t)$ as 
\begin{equation}
	s_{k,i}(t)= \sum_{j=0}^i \binom{k}{j} \dfrac{p_{k-j}^{(i-j)}(t)}{(i-j)!}.
\end{equation}

 Lower triangular matrices can be multiplied using their representation as differential operators given in Theorem \ref{diffOper}. 
\begin{thm}\label{mulForm}
 Let $A$ be a matrix of index zero and let $\{u_k(t)\}_{k\ge 0} $ be its associated sequence of polynomials. Let
	$$B=\sum_{i \ge 0} \dfrac{D^i}{i!} q_i(X),$$
where $q_i(t)$ is a polynomial of degree at most equal to $i$, for $i \ge 0$. Then we have
	\begin{equation}\label{eq:multAB}
A B= \sum_{k\ge 0} \dfrac{D^k}{k!} r_k(X),
	\end{equation}
	where
	$$r_k(t)= \sum_{i=0}^k s_{k,i}(t) q_i(t), \qquad k\ge 0,$$
	 and $s_{k,i}(t)$ is defined in \ref{eq:ski}.
\end{thm}

{\it Proof.} By \ref{eq:Aski} and the definition of $r_k(t)$ we have
\begin{eqnarray*}
	A B &=& \sum_{i\ge 0} A \dfrac{D^i}{i!} q_i(X) \cr 
	  &=& \sum_{i \ge 0} \sum_{k\ge 0} \dfrac{D^k}{k!} s_{k,i}(X) q_i(X) \cr
	  &=& \sum_{k \ge 0} \dfrac{D^k}{k!} \sum_{i=0}^k s_{k,i}(X) q_i(X) \cr
	  &=& \sum_{k\ge 0} \dfrac{D^k}{k!} r_k(X). 
	 \end{eqnarray*} \eop

	Let $A=\sum_{k\ge 0} \frac{D^k}{k!} p_k(X)$, where $p_k(t)$ is a polynomial of degree at most $k$, for $k \ge 0$. Then we have
\begin{equation}
	\dfrac{D^n}{n!} A = \sum_{k\ge n} \binom{k}{n} \dfrac{D^k}{k!} p_{k-n}(X), \qquad n \ge 0.
\end{equation}
From the Pincherle differentiation formulas \ref{eq:PinchDj} and \ref{eq:PinchXj} we obtain the Pincherle derivatives of $A$
\begin{equation}\label{eq:PinchXA}
X A - A X = \sum_{k\ge 0} \dfrac{D^k}{k!} p_{k+1}(X),
\end{equation}
and
\begin{equation}\label{eq:PinchAD}
A D - D A = \sum_{k\ge 1} \dfrac{D^k}{k!} p_k^\prime(X).
\end{equation}
Equation \ref{eq:XkDj} gives us
\begin{equation}\label{eq:XnA}
	X^n A= \sum_{k\ge 0} \dfrac{D^k}{k!} \sum_{j=0}^n \binom{n}{j} X^{n-j} p_{k+j}(X), \qquad n \ge 0.
\end{equation}

\subsection{Monomiality operators of general polynomial sequences}
Let $A$ be an invertible matrix of index zero, that is, an element of the group $\cG$. We obtain next some properties of the matrices $M=A^{-1} X A$ and $P=A^{-1} D A$. We will see how $M$ and $P$ are related with the Pincherle derivatives of $A$.

Let $\{u_k(t)\}_{k\ge 0}$ be the polynomial sequence associated with $A$ and let $p_k(t)$ for $k \ge 0$ be the polynomials that satisfy
$$A= \sum_{k\ge 0} \dfrac{D^k}{k!} p_k(X).$$

For a general element $A$ of $\cG$ the degree of $p_k(t)$ is at most $k$. If $A$ is a monic matrix in $\cG$, or the diagonal of $A$ of index zero is a multiple of the identity matrix, then, from the first equation in \ref{eq:binomTrans}, we can see that the degree of $p_k(t)$ is less that $k$, for $k \ge 1.$

Since $X \hX =I$ if $N=A^{-1} \hX A$ then we have $M N= I$, that is, $M$ has a right inverse. 

Let $J$ be the diagonal matrix of index -1 with entries $J_{k,k+1}=1/(k+1)$, for $k\ge 0$. Then we have $J D= I$. We call $J$ the {\em integration operator}.
Therefore $Q=A^{-1} J A$ is a left inverse of $P=A^{-1} D A$. 

Since $A^{-1} (X A - AX)= M-X $, from \ref{eq:PinchXA} we obtain
\begin{equation}
	M = X + A^{-1} \sum_{k\ge 0} \dfrac{D^k}{k!} p_{k+1}(X).
\end{equation}
Note that if $A$ is monic then the degree of $p_{k+1}(t)$ is at most equal to $k$. Therefore the sum in the previous equation is a matrix of index zero, and consequently, $M$ equals $X$ plus a matrix of index zero.

In an analogous way we obtain 
\begin{equation}
M=X - (X A^{-1} - A^{-1} X) A.
\end{equation}

Since $ A^{-1} ( A D - D A)= D - P$, from \ref{eq:PinchAD} we obtain
\begin{equation}
	P= D - A^{-1} \sum_{k\ge 1} \dfrac{D^k}{k!} p_k^\prime(X),
\end{equation}
and therefore, if $A$ is monic then $P$ equals $ D$ minus a matrix of index zero. 

The matrices analogous to the matrices $M$ and $P$ that we obtain by interchanging $A$ and $A^{-1}$ have some interesting properties. They are the monomiality matrices of the polynomial sequence associated with $A^{-1}$. 

\begin{thm}
Let $A$ be a monic element of $\cG$ and $\{u_k(t)\}_{k\ge 0}$ be polynomial sequence associated with $A$. Define $L= A X A^{-1}$.

For $m \ge 0$ define the linear functionals $\phi_m$ on the space of polynomials as follows:
\begin{equation}
	\phi_m (w(t))= w(L)_{0,m}, \qquad w(t) \in \cP.
\end{equation}
Then we have the biorthogonality relation
\begin{equation}
	\phi_m(u_n(t))=u_n(L)_{0,m}=\delta_{m,n}, \qquad m,n \ge 0,
\end{equation}
and therefore $\{\phi_m: m \ge 0\}$ is the dual basis of $\{u_n(t)\}_{n\ge 0}$.
\end{thm}
For a proof of this result see \cite[Cor. 3.3]{Lin}. A different construction of the dual basis was obtained in \cite{BenCh2}.

If $A$ is a monic element of $\cG$ then the matrix $L$ equals $X$ plus a lower triangular matrix. We will see in Section \ref{Appell} that when $A$ is an Appell matrix then $L+M= 2 X$.

If the sequence $\{u_n(t)\}_{n\ge 0}$ is orthogonal then $L$ is tridiagonal and it is called the Jacobi matrix of the orthogonal polynomial sequence.

Define $Q=A D A^{-1}$. Then we have $LQ - Q L =I$ 

\subsection{Differential operator representation of matrices of negative index}

Let $m$ be a positive integer and let $R$ be a matrix of order $-m$. Then we can write $R=A+B$, where $A$ is the lower triangular part of $R$ and $B=[b_{i,j}]$ is a banded matrix whose nonzero entries are in the diagonals of index in $\{-m, -m+1,\ldots, -1\}$. The diagonals of $B$ can be written in terms of the matrices $\hX$ and $\hD$ as follows. Define the numbers
\begin{equation}
e_{k,j}= \sum_{i=0}^{k-1} \binom{k}{i+j} (-1)^{k-i-j} b_{i,i+j}, \qquad 1 \le j \le m,\ k\ge j,
\end{equation}
and the polynomials
\begin{equation}
p_k(t)= \sum_{j=1}^k e_{k,j} t^{k-j}, \qquad k\ge 0.
\end{equation}
Then we have
\begin{equation}
B=\sum_{k=1}^m p_k(\hX) \dfrac{\hD^k}{k!}.
\end{equation}

\section{Monomiality and Appell polynomial sequences}\label{Appell} 

Let $A=[a_{k,j}]$ be a monic element of the group $\cG$ and let $\{u_k(t)\}_{k\ge 0}$ be the polynomial sequence associated with the rows of $A$. Then $A$ is invertible and $u_k(t)$ is a monic polynomial of degree $k$ for $k \ge 0$. In Section \ref{algcL} we saw that $A \rightarrow X A$ corresponds to the map $u_k(t)\rightarrow u_{k+1}(t)$, and also that $A \rightarrow D A$ corresponds to $u_k(t)\rightarrow k u_{k-1}(t)$. Now we use the simple Proposition \ref{basic}. Let $M=A^{-1} X A $ and $P=A^{-1} D A$. Then it is clear that $A M = X A$ and $A P = D A$. Since $\{u_k(t)\}_{k\ge 0}$ is a basis for the space of polynomials it is clear that multiplication of $A$ on the right by the matrices $M$ and $P$ defines linear operators on the space $\cP$ of polynomials in $t$. If we use $\hM$ and $\hP$ to denote the corresponding linear operators then we have
\begin{equation}\label{eq:mono}
	\hM u_k(t)= u_{k+1}(t), \quad \hbox{and} \quad \hP u_k(t)=k u_{k-1}(t), \qquad k \ge 0. 
\end{equation} 

Equation \ref{eq:mono} means that the polynomial sequence $\{u_k(t)\}_{k\ge 0}$ has the {\em monomiality property}.
 Therefore we have shown that every monic polynomial sequence has the monomiality property. 
 It is easy to see that, with a suitable re-scaling, we can show that every polynomial sequence has the monomiality property. See \cite{BenCh} for a different proof of this result. 
 
Let $\cF$ be the set of all formal power series of the form $g(t)=\sum_{k\ge 0} g_k t^k$, where the coefficients $g_k$ are complex numbers. The set $\cF$, with the usual addition and multiplication of power series, is a commutative algebra over the complex numbers.

 Let $\cF_0$ be the subset of $\cF$ of all the series $g(t)=\sum_{k=0}^\infty g_k t^k $, with $g_0\neq 0$. It is well-known that $\cF_0$ is a commutative group under the usual multiplication of formal power series.

The map $g(t) \rightarrow g(\hX)=\sum_{k\ge 0} g_k \hX^k$ is an algebra homomorphism from $\cF$ into $\cL$. The matrices $g(\hX)$ are Toeplitz matrices of index zero and they are invertible if $g_0\ne 0$. The image of $\cF_0$ under such map is the group of invertible Toeplitz matrices of index zero, and it is a subgroup of $\cG$. See \cite{Henr}.

Let us recall that $D= F \hX F^{-1}$, where $F=\hbox{Diag}(0!,1!,2!,\ldots)$. Then for every $g(t)$ in $\cF_0$ we have $F g(\hX) F^{-1}= g(D)$. Therefore the set of all matrices $g(D)$ with $g(t)$ in $ \cF_0$ is isomorphic to the group of invertible Toeplitz matrices of index zero.

 Let $f(t)=\sum_{k\ge 0} f_k \frac{t^k}{k!}$ be an element of the group $\cF_0$. We define the matrix 
 \begin{equation}\label{eq:f(D)}
	 f(D)= \sum_{k=0}^\infty f_k \dfrac{D^k}{k!}.
	\end{equation}
The map $f(t) \rightarrow f(D)$ is a linear homomorphism from the group $\cF_0$ into the group $\cG$.
Therefore we have $(f(D))^{-1}= (1/f)(D).$

Let $\{u_k(t)\}_{k\ge 0}$ be the polynomial sequence associated with the rows of $f(D)$.
 Then it is easy to see that 
 \begin{equation}\label{eq:binompoly}
	 u_k(t)=\sum_{j=0}^k \binom{k}{j} f_{k-j} t^j=f(\partial_t) t^k, \qquad k\ge 0,
	 \end{equation}
where $\partial_t$ denotes differentiation with respect to $t$.
From \ref{eq:binompoly} it follows that $u_k^\prime(t)= k u_{k-1}(t)$ for $k \ge 0$.

We also have the generating function
\begin{equation}\label{eq:genFct}
f(z) e^{tz}= \sum_{k=0}^\infty u_k(t) \dfrac{z^k}{k!}.
\end{equation}
Therefore $\{u_k(t)\}_{k\ge 0}$ is an {\em Appell polynomial sequence}. 

Since $f(D)$ commutes with $D$ it is clear that $P=(1/f)(D) D f(D)= D$.

Let $M=(1/f)(D)\, X \, f(D)$. Then we have $f(D) \, M= X \, f(D).$

\begin{thm}\label{logderiv}
	Let $f(t)=\sum_{k\ge 0} f_k \frac{t^k}{k!}$, with $f_0\ne 0$, and let
	 $h(t)=\frac{f^\prime(t)}{f(t)}$. Then we have  
	 $$M =\left(\dfrac{1}{f}\right)(D) X f(D) = X + h(D).$$ 
\end{thm}
{\it Proof.} 
 By the Pincherle differentiation formula $XD^k - D^k X= k D^{k-1}$, for $k\ge 1$ we get  
 $X f(D)-f(D) X= f^\prime(D)$ and then 
 $$(1/f)(D) f^\prime(D)=(1/f)(D) X f(D) - X= M - X,$$
and we obtain $M = h(D)+ X$. \eop

The previous Theorem shows that the linear operator $M$ that acts on the polynomials $u_k(t)$ has a representation as multiplication by $t$ plus a differential operator of infinite order, or of finite order, when $h(t)$ is a polynomial.  
 
Let $h(t)$ be as in Theorem \ref{logderiv} and let $\tilde{h}(t) \in \cF$ be an anti-derivative of $h(t)$, that is, $(\tilde{h})^\prime(t)=h(t)$. Then we have  $f(t) = c e^{\tilde{h}(t)}$ for some nonzero constant $c$ and therefore we obtain $f(D)= c e^{\tilde{h}(D)}$. This result can be used to construct examples of matrices $f(D)$ for which the operator $M$ is equal to $X$ plus a polynomial in $D$. If $h(t)$ is a polynomial of degree $m-1$ then $\tilde{h}(t)$ is a polynomial of degree $m$, $f(D)=f_0 e^{\tilde{h}(D)}$, and $M=X +h(D)$.
Let $\tilde{h}(t)=\sum_{k=1}^{m} y_k t^k$. Then
\begin{equation}\label{eq:severaly}
u_n(t)=\exp\left(\sum_{k=1}^{m} y_k \partial_t^k \right) t^n, \qquad n \ge 0,
\end{equation}
and this equation shows that $u_n$ is a polynomial function of $y_1,y_2,y_3,\ldots,y_m$ and $t$. A simple case is obtained with $\tilde{h}(t)= y_1 t + y_2 t^2$. In this case the polynomials $u_n(t)$ are generalized Hermite polynomials with two parameters. See \cite{Chag}, \cite{DatLic2022} and \cite{DLS2023} for some related results about generalized Hermite polynomials.

We will find next the matrices $f(D)$ for which $\{u_k(t)\}_{k\ge 0}$ is an orthogonal polynomial sequence. Let $f(t)$ be a series in $\cF_0$ and define $L=f(D) X (1/f)(D)$. Since the logarithmic derivative of $1/f$ is equal to $-h(t)$, where $h(t)=f^\prime(t)/f(t)$, by Theorem \ref{logderiv} we get $ L=X -h(D)$. From the equation $L f(D)= f(D) X$ we obtain the recurrence relation
\begin{equation}\label{eq:recrel}
	t u_n(t)= u_{n+1}(t) + \sum_{j=0}^n L_{n,j} u_j(t), \qquad n \ge 0.
\end{equation}

If $h(t)=h_0 + h_1 t$, with $h_1 \neq 0$, then $L= X -h_0 I - h_1 D$ is tridiagonal with nonzero entries in the diagonal of index one. In this case the recurrence relation \ref{eq:recrel} becomes the three-term recurrence relation
\begin{equation}\label{eq:trrr}
	t u_n(t) +n h_1 u_{n-1}(t) + h_0 u_n(t)-u_{n+1}(t)=0, \qquad n\ge 1.
\end{equation}
Therefore the polynomial sequence $\{u_k(t)\}_{k\ge 0}$ determined by $f(D)$ is orthogonal. If $h(t)$ is a polynomial of degree one then, in terms of the coefficients of $f(t)$ it is
\begin{equation}\label{eq:Orthologder}
h(t)=\dfrac{f_1}{f_0} + \dfrac{f_0 f_2 - f_1^2}{f_0^2} t,
\end{equation}
and therefore the Jacobi matrix $L$ is
\begin{equation}\label{eq:JacoD}
	L= X-\dfrac{f_1}{f_0} I - \dfrac{f_0 f_2 - f_1^2}{f_0^2} D,
\end{equation}
and the series $f(t)$ is
\begin{equation}\label{eq:fDortho}
	f(t)= f_0 \exp\left(\dfrac{f_1}{f_0}t + \dfrac{f_0 f_2 - f_1^2}{2 f_0^2} t^2\right)
\end{equation}

In this case the orthogonal polynomial sequence $\{u_k(t)\}_{k\ge 0}$ is a simple generalization of the Hermite polynomial sequence, since we have the parameters $f_1$ and $f_2$. The coefficient $f_0$ can be taken equal to one.

\section{Modified composition matrices and polynomial sequences of binomial type}

Let $f(z)=\sum_{k\ge 1} f_k z^k$, with $f_1=1$. Then $f^n(z)= z^n$ plus terms with higher powers of $z$. Let $C_f$ be the lower triangular matrix that has in the $n$-th column the coefficients of $f^n(z)$, that is, $(C_f)_{k,n} $ equals the coefficient of $z^k$ in the series $f^n(z)$. Let $g(z)=\sum_{k\ge 0} g_k z^k$. Then
$$C_f [g_0, g_1, g_2,\ldots ]\tr= [h_0,h_1, h_2, \ldots]\tr, $$
where 
$$h(z)= \sum_{k\ge 0} h_k z^k = g(f(z))=\sum_{k\ge 0} g_k f^k(z).$$
The matrix $C_f$ is called the composition matrix of $f(z)$, it is an invertible lower triangular matrix and its inverse is the composition matrix of a series $\tilde{f}(z)$ that satisfies $f(\tilde{f}(z))=z=\tilde{f}(f(z))$. The composition matrices form a non-commutative subgroup of $\cG$. See \cite{Henr} and \cite{Dual}.

We introduce next a set of modified composition matrices.
Let $f(z)=\sum_{k\ge 1} f_k \frac{z^k}{k!}$, with $f_1=1$,
and let
$$G(z,t)= \sum_{k\ge 0} f^k(z) \dfrac{t^k}{k!}=\exp( t f(z)).$$
The polynomial sequence $\{u_k(t)\}_{k\ge 0}$ is defined by the generating function:
\begin{equation}\label{Gzt}
G(z,t)= \sum_{k\ge 0} u_k(t) \dfrac{z^k}{k!}.
\end{equation}
Since $G(z,t)G(z,x)=G(z,t+x)$ it is easy to see that
$$u_n(t+x)= \sum_{k=0}^n \binom{n}{k} u_k(t) u_{n-k}(x), \qquad n \ge 0.$$
This means that $\{u_k(t)\}_{k\ge 0}$ is a {\em polynomial sequence of binomial type}. Such sequences have been extensively studied. See \cite{Roman} and \cite{Rota}. 

Let $B$ be the lower triangular matrix associated with the sequence $\{u_k(t)\}_{k\ge 0}$. Then $B$ is a monic element of the group $\cG$. It is easy to verify that
$$ B= F C_f F^{-1}, $$
where $C_f$ is the composition matrix of $f(z)=\sum_{k\ge 0} \frac{f_k}{k!} z^k$ and $F=\hbox{Diag}(0!,1!,2!,\ldots)$. Since the set of composition matrices is a subgroup of $\cG$ so is the set of all matrices of the form $F C_f F^{-1}$, where $C_f$ is a composition matrix. 

 It is clear that $B^{-1}=F C_{\tilde{f}} F^{-1}$, where $\tilde{f}$ is the inverse under composition of $f$. Therefore, since $F^{-1} X F=\hD$ we have 
 $$M= B^{-1} X B= F C_{\tilde{f}} F^{-1} X F C_f F^{-1}= F C_{\tilde{f}} \hD C_f F^{-1}.$$
 Let us compute the $n$-th column of $M$ using the properties of the factors in the previous equation. We identify a power series in $z$ with the infinite column vector of its coefficients. 
 Let $e_n$ be the infinite column vector that has its $n$-th coordinate equal to 1 and all the other coordinates equal to zero. If we identify $e_n$ with the series $z^n$ then $F^{-1} z^n= \frac{z^n}{n!}$ and $C_f F^{-1} z^n= \frac{f^n(z)}{n!}$. 

 Let us recall that multiplication of a column vector that corresponds to a series $g(z)$ on the left by $\hD$ gives the vector that corresponds to $g^\prime(z)$. Then we have 
 $$\hD \frac{f^n(z)}{n!}= f^\prime(z) \dfrac{f^{n-1}(z)}{(n-1)!},$$ 
 and multiplication on the left by $C_{\tilde{f}}$ yields
 $$C_{\tilde{f}} \left( f^\prime(z) \dfrac{f^{n-1}(z)}{(n-1)!}\right)= f^\prime(\tilde{f}(z)) \dfrac{z^{n-1}}{(n-1)!}.$$
 Since $f(\tilde{f}(z))=z$ we have $ f^\prime(\tilde{f}(z))= 1/{\tilde{f}^\prime(z)}.$
Let 
$$ \dfrac{1}{\tilde{f}^\prime(z)}= \sum_{k\ge 0} a_k \dfrac{z^k}{k!}.$$

 Finally, the $n$-th column of $M$ corresponds to the series
 $$M z^n= F \left( \sum_{k\ge 0} a_k \dfrac{z^{k+n-1}}{k! (n-1)!}\right)= \sum_{k\ge 0} \binom{k+n-1}{k} a_k z^{k+n-1},$$
 and therefore the entries in the diagonal of index $k-1$ of $M$ are 
 $$M_{n+k-1,n} = \binom{n+k-1}{k} a_k,$$
 and, since the $0$-th column of $M$ equals zero, this implies that
\begin{equation}\label{eq:MtildefDX}
	 M=\left( \sum_{k\ge 0} a_k \dfrac{D^k}{k!}\right) X=\left( \dfrac{1}{\tilde{f}^\prime}\right)(D)\, X.
	 \end{equation}
Note that $M$ determines a differential operator of infinite order with coefficients that are polynomials of degree one.

Now let $P=B^{-1} D B$. Then
$P= F C_{\tilde{f}} F^{-1} D F C_f F^{-1}.$
Since $ F^{-1} D F= \hX$ the $n$-th column of $P$ is given by
$$Pz^n= F C_{\tilde{f}} \hX \dfrac{f^n(z)}{n!} = F C_{\tilde{f}}\left( \dfrac{z f^n(z)}{n!}\right),$$
and then
$$Pz^n= F \left( \tilde{f}(z) \dfrac{z^n}{n!}\right).$$
Let $\tilde{f}(z)= \sum_{k\ge 1} g_k \frac{z^k}{k!}$. Then we have
$$Pz^n= F \sum_{k\ge 1} g_k \dfrac{z^{k+n}}{k! n!}=\sum_{k\ge 1} g_k \binom{k+n}{n} z^{k+n}.$$
This equation shows that, for $k\ge 1$, the entries of the diagonal of index $k$ of $P$ are given by
$$P_{n+k,n}= g_k \binom{n+k}{k}, \qquad k \ge 1, n\ge 0.$$
Therefore we have
\begin{equation}\label{eq:PtildefD}
P=\sum_{k\ge 1} g_k \dfrac{D^k}{k!}= \tilde{f}(D).
\end{equation}
This equation and \ref{eq:MtildefDX} give us immediately
$$ PM =\left( \dfrac{\tilde{f}}{\tilde{f}^\prime}\right)(D)\, X.$$

\section{Sheffer polynomial sequences}
In this Section we consider matrices associated with polynomial sequences called Sheffer sequences \cite{Roman}. The Sheffer polynomial sequences have been extensively studied using diverse approaches. See \cite{DMSri07}, \cite{Penson}, and \cite{Rota}. 

Let $g(t)=\sum_{k\ge 0} g_k \frac{t^k}{k!}$ with $g_0 \ne 0$ and let $A=g(D)$. Let $f(z)=\sum_{k\ge 1} f_k \frac{z^k}{k!}$, let $C_f$ be the composition matrix of $f$, and let $B=F C_f F^{-1}$. 

The matrix $S=BA$ is called Sheffer matrix \cite{Costa}. The set of all Sheffer matrices is a group and it is isomorphic to the group of the Riordan matrices, which are of the form $C_f g(\hX)$. Using the results in the previous Sections we can find easily the matrices $M_S=S^{-1} X S$ and $P_S=S^{-1} D S$, expressed as matrices of differential operators.

Let $\tilde{f}$ be the inverse under composition of $f$. By \ref{eq:MtildefDX} 
we have
$$M_S=S^{-1} X S= A^{-1} B^{-1} X B A = \left(\dfrac{1}{g}\right)(D) \left(\dfrac{1}{\tilde{f}^\prime}\right)(D) X g(D).$$
From the Pincherle differentiation formula \ref{eq:PinchDj} we get $X g(D)-g(D) X = g^\prime(D)$ and therefore 
\begin{eqnarray}
	M_S &=& \left(\dfrac{1}{g \tilde{f}^\prime }\right)(D) \left( g^\prime (D) + g(D) X \right)\cr
	 &=& \left( \dfrac{g^\prime}{g \tilde{f}^\prime}\right)(D) + \left(\dfrac{1}{\tilde{f}^\prime}\right)(D) X \cr
	 &=& \left(\dfrac{1}{\tilde{f}^\prime}\right)(D) \left( \left( \dfrac{g^\prime}{g }\right)(D)+X \right).
	 \end{eqnarray}

By \ref{eq:PtildefD} we have 
$$P_S= S^{-1} D S = A^{-1} B^{-1} D B A = \left( \dfrac{1}{g}\right)(D) \tilde{f}(D) g(D) = \tilde{f}(D).$$

Therefore
$$P_S M_S =\left(\dfrac{\tilde{f}}{\tilde{f}^\prime}\right)(D) \left(\left( \dfrac{g^\prime}{g }\right)(D) + X\right).$$

Since the group of Riordan matrices is isomorphic to the group of Sheffer matrices we can easily study the monomiality properties of the polynomial sequences associated with Riordan matrices. 

\section{Some orthogonal polynomial sequences}
In this Section we consider a set of orthogonal polynomial sequences that are associated with invertible lower triangular matrices $A$ for which the matrices $M=A^{-1} X A$ and $P=A^{-1} D A$ have a relatively simple structure. 

Let us define polynomials $f(t)=f_0 + f_1 t + f_2 t^2$ and $g(t)= g_0 + g_1 t$, where $g_1 \ne 0$, and the matrix 
\begin{equation}\label{eq:DiffOper}
	B= D g(X) + \dfrac{D^2}{2} f(X).
\end{equation}
Note that $B $ is a banded matrix of index zero. Define
\begin{equation}
H= g_1 D X + f_2 \dfrac{D^2}{2} X^2.
\end{equation}
The matrix $H$ is diagonal of index zero and coincides with the diagonal of index zero of $B$.
Let $A$ be a monic matrix of index zero that satisfies the equation $A B = H A$.
Then $L=A X A^{-1}$ is a tridiagonal matrix of index -1 and the polynomial sequence $\{u_k(t)\}_{k\ge 0}$ associated with $A$ is orthogonal with respect to a linear functional on the space of polynomials. See \cite{Opm} and \cite{Uni}. 

The collection of all the polynomial sequences associated with matrices $A$ constructed as described above includes almost all the hypergeometric orthogonal polynomial sequences. See \cite{Hyp}. We consider next the subset of the polynomial sequences obtained as described above when the coefficient $f_2$ is equal to zero.

If $f_2=0$ the tridiagonal Jacobi matrix $L=A X A^{-1}$ is
\begin{equation}\label{eq: Lsimple}
L= X -\dfrac{g_0}{g_1} I -\dfrac{f_1}{g_1} D X +\dfrac{f_1 g_0 -f_0 g_1}{2 g_1^2} D + \dfrac{f_1^2}{4 g_1^2} D^2 X.
\end{equation}
Here we have the coefficients $f_0,f_1, g_0,g_1$ as parameters. Giving appropriate values to the parameters we can obtain the Jacobi matrices of the Meixner-Pollaczek, Meixner, Krawtchouk, Laguerre, and Hermite orthogonal polynomial families. See \cite{Hyp} and \cite{SriBCh}. 

When $f_2=0$ the matrix $M=A^{-1} X A$ is tetradiagonal of index -1 and 
\begin{equation}\label{Mortho}
	M= w_0(X) + D w_1(X) + \frac{1}{2} D^2 w_2(X),
\end{equation}
where
\begin{eqnarray}
	w_0(t)&=& \dfrac{g_0}{g_1} + t, \cr
	w_1(t)&=& \dfrac{f_0 g_1 + f_1 g_0}{2 g_1^2} +\dfrac{f_1}{g_1} t,\cr
	w_2(t)&=& \dfrac{f_0 f_1}{2 g_1^2} + \dfrac{f_1^2}{2 g_1^2} t.
\end{eqnarray}
The linear operator on the space of polynomials determined by $M$ is multiplication by the variable $t$ plus a differential operator of order two with polynomial coefficients of order one.

The matrix $P=A^{-1} D A$ in this case is 
\begin{equation}
	P=\sum_{k\ge 0} \left(\dfrac{-f_1}{2 g_1}\right)^k D^{k+1}= D \exp\left(\dfrac{-f_1}{2 g_1} D \right).
\end{equation}
Note that the linear operator on the space of polynomials determined by $P$ is differentiation followed by a translation operator, and such operations commute. 

It is well-known that a change in the value of the parameter $g_1$ corresponds to a linear change of variable in the orthogonal polynomials $u_k(t)$. Two orthogonal polynomial sequences that differ by a linear change of variable are usually considered as equivalent sequences.  

The case of the generalized Laguerre polynomials is discussed in the next section, where the corresponding matrix is obtained as a series of powers of a generalized differentiation operator.

In the cases with $f_2\ne 0$ the matrices $M$ and $P$ are more complicated, but can be obtained with the methods presented above.

\section{Series of generalized derivatives}

In this Section we consider a class of polynomial sequences determined by the rows of invertible matrices that are power series of a generalized derivative $D_c$. We will show that the results of Section \ref{Appell} can be easily generalized and can be applied to a larger class of polynomial sequences. The generalized derivatives that we consider here are usually called Ward derivatives. 

Let $c_1,c_1,c_3,\ldots $ be a sequence on nonzero complex numbers. Define the function $f_c$ on the set on non-negative integers by $f_c(n)=c_1 c_2 \cdots c_n$, for $n \ge 1$, and $f_c(0)=1$. We call $f_c$ the $c$-factorial function.
The $c$-binomial coefficients are defined by 
$$cbin(n,k)=\dfrac{f_c(n)}{f_c(k) f_c(n-k)}, \qquad n \ge k \ge 0.$$

Define the diagonal matrix of index zero
\begin{equation}\label{eq:Fc}
	 F_c=\hbox{Diag}\left(\dfrac{f_c(0)}{0!}, \dfrac{f_c(1)}{1!}, \dfrac{f_c(2)}{2!}, \ldots \right).
\end{equation}

The $c$-derivative matrix $D_c$ is defined by $D_c=F_c D F_c^{-1}$. It is a diagonal matrix of index one and its entries are $(D_c)_{k+1,k}=c_{k+1}$, for $k \ge 0$. 

The $c$-shift matrix $X_c$ is defined by $X_c=F_c X F_c^{-1}$. It is a diagonal matrix of order -1 and its entries are $(X_c)_{k,k+1}= \frac{k+1}{c_{k+1}}$, for $ k \ge 0$. 
Since $XD-DX=I$, it is clear from the definitions of $Dc$ and $X_c$ that $X_c D_c - D_c X_c=I$.

Let 
$$ g(t)= \sum_{k=0}^\infty g_k \dfrac{t^k}{f_c(k)}, $$
be such that $g_0 \ne 0$, that is, $g(t)$ is in the group $\cF_0$.
Then the matrix $g(D_c)$ is a well-defined element of the group $\cG$. 
Let $\{u_k(t)\}_{k \ge 0}$ be the polynomial sequence determined by the rows of $g(D_c)$. Then we have
\begin{equation}\label{eq: cbinpoly}
	u_n(t) = \sum_{k=0}^n cbin(n,k) g_{n-k} t^k, \qquad n \ge 0. 
\end{equation}

The $c$-exponential series is 
\begin{equation}\label{eq:cexp}
	e_c(t) = \sum_{k=0}^\infty \dfrac{t^k}{f_c(k)}.
\end{equation}

It is easy to verify that the sequence $\{u_k(t)\}_{k\ge 0}$ has the following generating function
\begin{equation}\label{eq:cgenfct}
	g(z) e_c(tz) =\sum_{k=0}^\infty u_k(t) \dfrac{z^k}{f_c(k)}.
\end{equation}
We say that $\{u_k(t)\}_{k\ge 0}$ is a sequence of generalized Appell type.

Since $g(D_c)$ commutes with $D_c$ we have $P=(1/g)(D_c) D_c g(D_c)=D_c.$

We define the matrix $M=(1/g)(D_c) X_c g(D_c)$. It is clear that $g(D_c) M = X_c g(D_c)$ and therefore the operator $M$ sends the polynomial $u_k(t)$, which corresponds to the $k$-th row of $g(D_c)$, to $\frac{k+1}{c_{k+1}} u_{k+1}(t)$, that corresponds to the $k$-th row of $X_c G(D_c)$.

Since we have the identity $X_c D_c-D_c X_c=I$, the proof of the following Theorem is analogous to the proof of Theorem \ref{logderiv}.

\begin{thm}\label{clogderiv}
Let $h(t)=\frac{g^\prime(t)}{g(t)}$. Then we have $M =X_c + h(D_c)$.
\end{thm}

As in the case of matrices that are power series in $D$, the matrices $g(D_c)$ for which the operator $M$ is equal to $X_c$ plus a polynomial in $D_c$ are those for which the logarithmic derivative of $g(t)$ is a polynomial. 

Let $L=g(D_c) X_c (1/g)(D_c)$. In this case we also have $L-X_c=- (M-X_c) $, since the logarithmic derivative of $1/g(t)$ is the negative of the logarithmic derivative of $g(t)$. Therefore $L$ is tridiagonal if $g^\prime(t)/g(t)$ is a polynomial of degree one. In such case we obtain
\begin{equation}\label{eq:cL}
	L= X_c -\dfrac{g_1}{g_0 c_1} I - \dfrac{2 c_1 g_0 g_2- c_2 g_1^2}{c_1^2 c_2 g_0^2 } D_c.
\end{equation}
This matrix satisfies the equation $L g(D_c) = g(D_c) X_c$, which looks like the equation that gives the three-term recurrence relation of an orthogonal polynomial sequence, but the linear operator on the rows of $g(D_c)$ determined by multiplication on the right by $X_c$ is, in the general case, not equal to the operator of multiplication by the variable $t$. That occurs only when $D_c=D$, that is, when $c_k=k $, for $k\ge 1$. Therefore, if $D_c \ne D$ then the tridiagonal matrix $L$ given in \ref{eq:cL} is not the Jacobi matrix of an orthogonal polynomial sequence. 

If the matrix $L$ is tridiagonal then the row polynomials $u_k(t)$ of $g(D_c)$ are related to an orthogonal polynomial sequence $\{v_k(t)\}_{k \ge 0}$ as  we show next. Let $f_k=(k!/f_c(k)) g_k$, for $k \ge 0$, and let 
$$f(D)=\sum_{k=0}^\infty \dfrac{f_k}{k!} D^k.$$
Then we have $g(D_c)=F_c f(D) F_c^{-1}$. Since the logarithmic derivative of $g(t)$ is a polynomial of degree one, so is the logarithmic derivative of $f(t)$ and by the results in the previous Section the row polynomials $v_k(t)$ of the matrix $f(D)$ are orthogonal. If $v_k(t)=\sum_{j=0}^k a_{k,j} t^j $ then, from $g(D_c)=F_c f(D) F_c^{-1}$ we obtain 
\begin{equation} \label{eq:ukvk}
u_k(t)=\sum_{j=0}^k \dfrac{f_c(k) j!}{f_c(j) k!} a_{k,j} t^j, \qquad k\ge 0.
\end{equation}
Therefore the polynomials $u_k(t)$ determined by the matrix $g(D_c)$ are obtained by a modification of a sequence of generalized Hermite polynomials.

Let us note that the results in this Section hold for any sequence $c_1,c_2,c_3,\ldots $ of nonzero numbers.
We consider next some examples that correspond to particular choices of the sequence $c_1,c_2,c_3,\ldots $.

\bigskip
{\bf Example 1.}
If we take $c_k=k^2$, for $k\ge 1$, then $(D_c)_{k+1,k}=(k+1)^2$, for $k \ge 0$, and $(X_c)_{k,k+1}=1/(k+1)$, for $k \ge 0$. Let us note that $X_c D=I$, that is $Xc$ is a left inverse of $D$. Note also that $D_c= D X D$.
In this case we have $f_c(n)=(1)(4)(9)\cdots (n^2)$, for $n\ge 1$.

Let $$g(t)= \sum_{k=0}^\infty g_k \dfrac{t^k}{f_c(k)},$$
where $g_0\ne 0$, be an element of the group $\cF_0$. Define the matrix  
 $$g(D_c)= \sum_{k=0}^\infty g_k \dfrac{D_c^k}{f_c(k)}, $$
and let $\{u_k(t)\}_{k\ge 0}$ be the corresponding polynomial sequence.

If $g(t)$ is such that its logarithmic derivative $h(t)$ is a polynomial of degree one then
$$h(t)= \dfrac{g_1}{g_0} + \dfrac{g_0 g_2 - 2 g_1^2}{2 g_0^2} t,$$
and therefore
$$M=X_c + \dfrac{g_1}{g_0} I + \dfrac{g_0 g_2 - 2 g_1^2}{2 g_0^2} D_c,$$
and the matrix $L$ is
$$L=X_c - \dfrac{g_1}{g_0} I - \dfrac{g_0 g_2 - 2 g_1^2}{2 g_0^2} D_c. $$
We also have
\begin{equation}\label{eq:expInth}
	g(t)=g_0 \exp\left(\dfrac{g_1}{g_0} t + \dfrac{g_0 g_2 - 2 g_1^2}{4 g_0^2} t^2\right).
\end{equation}

The polynomial sequence $\{u_k(t)\}_{k\ge 0}$ obtained in this case satisfies the relation
\begin{equation}\label{eq:recrelD2}
	\dfrac{2 g_1^2-g_0 g_2}{2 g_0^2} k^2 u_{k-1}(t) -\frac{g_1}{g_0} u_k(t) +\dfrac{1}{k+1} u_{k+1}(t)= J_c u_k(t), \qquad k \ge 1,
\end{equation}
where $J_c$ denotes the indefinite integration operator, with integration constant equal to zero, which is the linear operator on the space of polynomials that corresponds to multiplication on the right by $X_c$.
 Notice that \ref{eq:recrelD2} does not have the form of the three-term recurrence relation of an orthogonal polynomial sequence. The polynomials $u_k(t)$ are related to orthogonal polynomials $v_k(t)$ of Hermite-type by \ref{eq:ukvk}. If we set $g_0=1, g_1=z$, and $ g_2=2 z^2+ 4 y$ then $$g(t)= \exp\left(z t + y t^2\right), $$ 
 the polynomials $u_k(t)$ are functions of $y, z$, and $t$, and $M= X_c + z I +y D_c$. We also have $u_0(t)=1$, and 
 $$ u_1(t)=t+z,\ u_2(t)=t^2+4zt+ 2 z^2+ 4 y, \ u_3(t)=t^3 +9 z t^2+(18 z^2 + 36 y)t +6z^3+36 y z.$$

\bigskip
{\bf Example 2.}
Let $c_k= a k+ y k^2$ for $k\ge 1$, where $a$ and $y$ are complex numbers and $a\ne 0$.  
Let $U$ be the diagonal matrix of index one with entries $U_{k+1,k}=(k+1)^2$, for $k\ge 0$. Then we have $D_c= a (D + y U)$. 
Notice that $D_c$ is diagonal of index one.

Let us define 
\begin{equation}\label{eq:GenLaguerre}
A= \exp(-D_c).
\end{equation}
Then $A$ is a monic element of $\cG$ and $A^{-1}=\exp(D_c)$. 
A simple computation yields
$$ L= A X A^{-1}=X + a(1+ y) I + D (a^2 y (1+ y) + 2 a y X) + a^2 y^2 D^2 X, $$
and
$$ M=A^{-1} X A = X - a(1+ y) I + D (a^2 y (1 +y) - 2 a y X) + a^2 y^2 D^2 X. $$ 
Note that in this case $L$ and $M$ are tridiagonal and hence the polynomial sequences associated with $A$ and $A^{-1}$ are orthogonal. Note also that $M$ is obtained from $L$ by substitution of $a$ with $-a$.

If $y=1/a$ then $L= (1+a) I + D ((1+a) I + 2 X) + D^2 X $ is the Jacobi matrix of the monic Laguerre polynomials. See \cite{Hyp}. Therefore the polynomial sequence $\{u_k(t)\}_{k\ge 0}$ associated with $A$ is a generalization of the Laguerre polynomials that has $a$ and $y$ as parameters. In \cite{SriBCh} some results about Laguerre polynomials are obtained using umbral and operational methods. See also \cite{RicciT}.

In this case the matrix $P=A^{-1} D A$ is
$$P=\sum_{k\ge 0} a^k y^k D^{k+1}, $$
and thus $P=f(D)$, where $f(t)=t/(1- a y t)$.

This example can be easily generalized by taking $c_k= \sum_{j=1}^m a_j k^j$, but the polynomial sequences obtained with such $c_k$ are not orthogonal if $m>2$. 

\bigskip
 {\bf Example 3.} 
 Now we consider the Jackson $q$-derivative $\cD_q$, defined by 
\begin{equation}\label{eq:Dq}
	\cD_q p(t)=\dfrac{p(q t)- p(t)}{q t - t}, \qquad p(t) \in \cP.
\end{equation}

 Let $q$ be a complex number that is not a root of unity, and define the basic numbers, or $q$-numbers, by
 \begin{equation}
	 [n]=\dfrac{q^n -1}{q-1}, \qquad n \in \Z.
	 \end{equation}
Let $c_k=[k]$, for $k \ge 1.$ Then the $c$-factorial function is $f_c(n)=[1][2][3]\cdots [n]$, which we denote by $[n]!$,
 the matrix $D_c$ is diagonal of index one with entries $(D_c)_{k+1,k}=[k+1]$, for $k \ge 0$, and $X_c$ is the diagonal matrix of index -1 with entries $(X_c)_{k,k+1}=(k+1)/[k+1]$, for $k \ge 0$. Let $\tilde{D}_c$ be the diagonal of index -1 with entries $(\tilde{D}_c)_{k,k+1}=1/[k+1]$, for $ k \ge 0$. Then $\tilde{D}_c D_c=I$ and $X_c=X D \tilde{D}_c=\tilde{D}_c D X$.

Let $f_0=1$,   
$$f(t)=\sum_{k=0}^\infty f_k \dfrac{t^n}{[n]!},$$
and
$$f(D_c) =\sum_{k=0}^\infty f_k \dfrac{D_c^n}{[n]!}.$$

Then the matrix $f(D_c)$ is lower triangular and invertible and its inverse is $(1/f)(D_c)$. In this case we have $P= D_c$ and $ M= X_c + h(D_c)$, where $h(t)$ is the logarithmic derivative of $f(t)$. 

A simple example is obtained when $h(t)$ is a polynomial of degree one. In this case we have  
$$h(t)= f_1+ \dfrac{2 f_2 -[2] f_1^2}{[2]} t, $$
and then
$$f(t)=\exp\left( f_1 t + \dfrac{2 f_2 -[2] f_1^2}{ 2 [2]} t^2\right),$$
and 
$$ M= X_c + f_1 I + \dfrac{2 f_2 -[2] f_1^2}{[2]} D_c.$$
The corresponding polynomials $u_k(t)$ have $f_1$ and $f_2$ as parameters and can be expressed as a transformation of generalized Hermite polynomials using \ref{eq:ukvk}.


\begin{thebibliography}{00} 
\bibitem{BenCh} Youss\`ef Ben Cheikh, Some results on quasi-monomiality, Appl. Math. Comput. 141 (2003) 63--76.
\bibitem{BenCh2} Youss\`ef Ben Cheikh, On obtaining dual sequences via quasi-monomiality, Georgian Math. J. 9 (2002) 413--422.
\bibitem{Chag} Hamza Chaggara, Operational rules and a generalized Hermite polynomials, J. Math. Anal. Appl. 332 (2007) 11--21.
\bibitem{Costa} Francesco A. Costabile, Modern Umbral Calculus, De Gruyter, Berlin/Boston, 2019.
\bibitem{Henr} Peter Henrici, Applied and computational Complex Analysis, Vol. I, J. Wiley, New York, 1974.
\bibitem{DCS97} G. Dattoli, C. Cesarano, M. Sacchetti, A note in the monomiality principle and generalized polynomials, J. Math. Anal. Appl. 227 (1997) 98--111.
\bibitem{DatLic2022} Giuseppe Dattoli, Silvia Licciardi, Monomiality and a new family of Hermite polynomials, arXiv:2205.11517v1 [math.CA].
\bibitem{DLS2023} Giuseppe Dattoli, Silvia Licciardi, Elio Sabia, An operational point of view to the theory of multi-variable/multi-index Hermite polynomials, arXiv:2310.19689v1 [math-ph].
\bibitem{DMSri07} G. Dattoli, M. Migliorati, H. M. Srivastava, Sheffer polynomials, monomiality principle, algebraic methods and the theory of classical polynomials, Mathematical Computer Modelling 45 (2007) 1033--1041.
\bibitem{Hyp} R. Koekoek, P. A. Lesky, R. F. Swarttouw, Hypergeometric Orthogonal Polynomials and Their $q$-Analogues, Springer Monographs in Mathematics, Springer-Verlag, Berlin, Heidelberg, 2010.
\bibitem{LiccDat} Silvia Licciardi, Giuseppe Dattoli, Guide to the Umbral Calculus: A different mathemetical language, World Scientific, Singapore, 2022.
\bibitem{Penson} K. A. Penson, P. Blasiak, G. Dattoli, G. H. E. Duchamp, A. Horzela, A. I. Solomon, Monomiality principle, Sheffer type polynomials and the normal ordering problem, in: The Eighth International School on Theoretical Physics “Symmetry and Structural Properties of Condensed Matter”, Myczkowce, Poland, August 2005, J. Phys. Conf. Ser. 30 (2006) 86--97.
\bibitem{RicciT} P. E. Ricci, I. Tavkhelidze, An introduction to operational techniques and special polynomials, J. Mathematical Sciences, 157 (2009) 161--189.
\bibitem{Roman} Steven Roman, The Umbral Calculus, Academic Press, 1984.
\bibitem{Rota} Gian-Carlo Rota, Finite Operator Calculus, Cambridge University Press, 2005.
\bibitem{SriBCh} H. M. Srivastava, Y. Ben Cheikh, Orthogonality of some polynomial sets via quasi-monomiality, Appl. Math. Comput. 141 (2003) 415--425.
\bibitem{Lin} L. Verde-Star, Linearization and connection coefficients of 	polynomial sequences: A matrix approach, Linear Algebra Appl. 672 (2023) 195--209.
\bibitem{Opm} L. Verde-Star, Characterization and construction of classical orthogonal polynomials using a matrix approach, Linear Algebra Appl. 438 (2013) 3635--3648.
\bibitem{PSGH} L. Verde-Star, Polynomial sequences generated by infinite Hessenberg matrices, Spec. Matrices 2017; 5: 64--72.
\bibitem{Uni} L. Verde-Star, A unified construction of all the hypergeometric and basic hypergeometric families of orthogonal polynomial sequences, Linear Algebra Appl. 627 (2021) 242--274.
\bibitem{Dual} L. Verde-Star, Dual operators and Lagrange inversion in several variables, Advances in Math. 58 (1985) 89--108.
\end{thebibliography}
\end{document}